\documentclass[psamsfonts]{amsart} 
\usepackage{amsmath, amssymb, url, algorithm, tikz, amsthm}
\usepackage[noend]{algorithmic}
\usepackage{subfig}
\usetikzlibrary{arrows,shapes}

\newcommand{\ai}{\mathfrak{a}}
\newcommand{\ah}{\hat{\alpha}}
\newcommand{\bi}{\mathfrak{b}}
\newcommand{\ch}{\mathrm{char}}

\newcommand{\D}{\mathcal{D}}
\newcommand{\dv}{\mathrm{div}}
\newcommand{\f}{\mathfrak{f}}
\newcommand{\g}{\mathfrak{g}}
\newcommand{\F}{\mathbb{F}}
\newcommand{\Fq}{\F_q}
\newcommand{\HH}{\mathcal{H}}
\newcommand{\I}{\mathcal{I}}
\newcommand{\idealp}{\mathrm{id}}
\newcommand{\ideal}[1]{\left\langle {#1} \right\rangle }

\newcommand{\J}{\mathcal{J}}

\newcommand{\kk}{\mathfrak{k}}
\newcommand{\N}{\mathbb{N}}
\newcommand{\OO}{\mathcal{O}}
\newcommand{\p}{\mathfrak{p}}
\newcommand{\PP}{\mathcal{P}}

\newcommand{\ph}{\hat{\psi}}
\newcommand{\ppmod}[1]{\mbox{\ $\left(\mathrm{mod}\ {#1}\right)$}}
\newcommand{\R}{\mathcal{R}}
\newcommand{\RR}{\mathbb{R}}
\newcommand{\st}{\mid}
\newcommand{\supp}{\mathrm{supp}}
\newcommand{\SC}{\mathcal{S}}
\newcommand{\tk}{\mathfrak{t}}
\newcommand{\wk}{\mathfrak{w}}
\newcommand{\Z}{\mathbb{Z}}
\newcommand{\reduce}{\mathop{\mathrm{red}}}

\newcommand{\WHILEOL}[2]{\STATE \textbf{while}\ #1\ \textbf{do}\ #2}
\newcommand{\IFOL}[2]{\STATE \textbf{if}\ #1\ \textbf{then}\ #2}

\newtheorem{theorem}{Theorem}

\newtheorem{proposition}{Proposition}

\theoremstyle{definition}

\begin{document}

\title[Computations in Cubic Function Fields of Unit Rank Two]{Class Number and Regulator
Computation in \\ Purely Cubic Function Fields of Unit Rank Two}

\author[F. Fontein]{Felix Fontein} \address{Felix Fontein, Department of Mathematics \&\ Statistics,
University of Calgary, 2500 University Drive NW, Calgary, Alberta, Canada T2N 1N4}
\email{fwfontei@ucalgary.ca}

\author[E. Landquist]{Eric Landquist} \address{Eric Landquist, Kutztown University, Department of
Mathematics, Kutztown, PA 19530, USA} \email{elandqui@kutztown.edu}

\author[R. Scheidler]{Renate Scheidler} \address{Renate Scheidler, Department of Mathematics \&\
Statistics, University of Calgary, 2500 University Drive NW, Calgary, Alberta, Canada T2N 1N4}
\email{rscheidl@ucalgary.ca}

\maketitle


\begin{abstract}
  We describe and give computational results of a procedure to compute the divisor class number and
  regulator of most purely cubic function fields of unit rank $2$.  Our implementation is an
  improvement to Pollard's Kangaroo method in infrastructures, using distribution results of class
  numbers as well as information on the congruence class of the divisor class number, and an
  adaptation that efficiently navigates these torus-shaped infrastructures.  Moreover, this is the
  first time that an efficient ``square-root'' algorithm has been applied to the infrastructure of a
  global field of unit rank $2$.  With the exception of certain function fields defined by Picard
  curves, our examples are the largest known divisor class numbers and regulators ever computed for
  a function field~of~genus~$3$.
\end{abstract}

\section{Introduction and Motivation}
\label{S:intro} 
One of the more difficult problems in arithmetic geometry is the computation of the divisor class
number of an algebraic curve over a finite field.  In this paper, we give results on the application
and optimization of a method of Scheidler and Stein \cite{ss07,ss08}, combined with modifications to
Pollard's Kangaroo algorithm \cite{pollardkangaroo}, to compute the exponent of the infrastructure
of a purely cubic function field with complete splitting at infinity (i.e., unit rank~$2$) over a
large base field. The regulator and the divisor class number are multiples of this exponent, and in
many cases all three numbers are the same, whence our algorithm computes the regulator and the
divisor class number in these cases as well. Our method greatly improves upon the method described
in \cite{lsy} to compute the regulator in this setting and is the first ever treatment of an
efficient ``square-root'' algorithm in a two-dimensional infrastructure of a global field.

An algorithm due to Stein and Williams \cite{sw} uses techniques of Lenstra \cite{len} and Schoof
\cite{schoof82} to compute the divisor class number and regulator of a real quadratic function field
in $O\left(q^{[(2g-1)/5] + \varepsilon(g)}\right)$ infrastructure operations, \footnote{Throughout
this paper $[r]$ will denote the nearest integer to $r \in \RR$.}  where $-1/4\leq \varepsilon(g)
\leq 1/2$.  This method was improved by Stein and Teske \cite{st02,st02b,st05}, who applied the
Kangaroo algorithm to compute the $29$-digit class number and regulator of a real quadratic function
field of genus~$3$.

The algorithm of \cite{sw} was generalized to cubic and arbitrary function fields in
\cite{ss07,ss08}, respectively, and implemented in purely cubic function fields of unit rank $0$ and
$1$ in \cite{ericphd}. In this paper, we provide an implementation and numerical examples for purely
cubic function fields of unit rank~$2$.  Our method is applied to compute divisor class numbers and
regulators of up to $31$ digits of function fields of genus $3$. With the exception of the
$55$-digit class numbers computed by Bauer, Teske, and Weng in \cite{btw,weng} for cubic function
function fields generated by Picard curves, our examples are the largest known class numbers and
regulators ever computed for a function field of genus at least~$3$ over a large base field.

The remainder of this paper is organized as follows. We first give an overview of cubic function
fields and their infrastructure. Then we outline our version of Pollard's Kangaroo method for
infrastructures, and explain how to recover the regulator and divisor class number in most
cases. Finally, we discuss details of our implementation and provide numerical results.

\section{Cubic Function Fields}
\label{S:background}

For a general introduction to function fields, we direct the reader to \cite{h,sti,rosen}. Explicit
details of purely cubic function fields and their arithmetic can be found in
\cite{ss00,s01,s04,b,ss07}. Let $\Fq$ be a finite field and $\Fq(x)$ the field of rational functions
in $x$ over $\Fq$.  Throughout this paper, we assume that $\ch(\Fq) \geq 5$. A {\em cubic function
field} is a separable extension $K/\Fq(x)$ of degree $3$; we denote by $g$ the genus of $K$. A
function field is \emph{purely} cubic if it is of the form $K = \F_q(x,y)$ where $y^3 = F$ for some
cube-free $F \in \F_q[x]$.

\subsection{Divisors and Ideals}
\label{S:divisors}

Let $\D$ denote the group of divisors of $K$ defined over $\Fq$, $\D_0$ the subgroup of divisors of
degree $0$ defined over $\Fq$, and $\PP$ the subgroup of principal divisors defined over $\Fq$. Then
the {\em $($degree~$0)$ divisor class group} of $K$ is the quotient group $\J = \D_0/\PP$ and its
order $h = |\J|$ is the {\em $($degree~$0)$ divisor class number} of $K$. Let $S$ be the set of
places of $K$ lying above the place at infinity of $\Fq(x)$, $\supp(D)$ the support of $D \in \D$,
$\D_0^S = \{D\in\D_0\st \supp(D) \subseteq S\}$, and $\PP^S = \PP\cap\D_0^S$.  Then the order $R$ of
the quotient group $\D_0^S/\PP^S$ is the {\em $(S$-$)$regulator} of $K$.  Finally, let $\D_S =
\{D\in\D\st \supp(D)\cap S = \emptyset\}$ and $\PP_S = \PP\cap\D_S$. Then every $D\in\D$ can be
uniquely written in the form $D = D_S + D^S$ with $D_S\in\D_S$ and $D^S\in\D^S$.

The {\em maximal order} of $K/\Fq(x)$ is the integral closure of $\Fq[x]$ in $K$ and is denoted
$\OO$.  Let $\I$ denote the group of non-zero fractional ideals of $\OO$ and $\HH$ the subgroup of
non-zero principal fractional ideals. The {\em ideal class group} of $K$ is the quotient group
$Cl(\OO) = \I/\HH$, and its order $h_{\OO} = |Cl(\OO)|$ is called the {\em ideal class number} of
$K$. Let $f$ be the greatest common divisor of the degrees of all the places in $S$. By Schmidt
\cite{schmidt} (see also Proposition~14.1 of \cite{rosen}) there is an exact sequence $$ (0)
\longrightarrow \D_0^S/\PP^S \longrightarrow \J \longrightarrow Cl(\OO) \longrightarrow \Z/f\Z
\longrightarrow (0) \enspace , $$ so that $fh = h_{\OO}R$.

There is a well-known isomorphism $\Phi: \D_S \to \I$ given by $D \mapsto \{\alpha \in K^* \st
\dv(\alpha)_S \geq -D\}$ with inverse $\f \mapsto -\sum_{\p\notin S} m_{\p}\p$, where $\p$ denotes
any finite place of $K$, $m_{\p} = \min\{v_{\p}(\alpha) \st 0\neq\alpha \in \f\}$, and $v_{\p}$ is
the normalized discrete valuation corresponding to $\p$. Moreover, $\Phi$ induces an isomorphism
from $\D_S/\PP_S$ to $Cl(\OO)$.  If $S$ contains an infinite place $\infty_0$ of degree $1$, then
$\Phi$ can be extended to an isomorphism $$\Psi: \big \{ \, D \in \D_0 \st v_{\p}(D) = 0 \mbox{ for
all } \p \in S \setminus \{ \infty_0\} \, \big \} \to \I$$ by $\Psi(D_S - \deg(D_S)\infty_0) =
\Phi(D_S)$, with the inverse given by $\Psi^{-1}(\f) = \Phi^{-1}(\f) +
\deg(N_{K/\Fq(x)}(\f))\infty_0$.

\subsection{Units}
\label{S:units}

By Proposition~14.1 and its Corollary~1 of \cite{rosen}, $\OO^* / \Fq^* \cong \PP^S$ is a free
abelian group of rank~$r = |S| - 1$. We write $S = \{\infty_0,\ldots,\infty_r\}$ to denote the
infinite places of $K$, with $v_i$ the normalized discrete valuation corresponding to $\infty_i$,
for $0\leq i \leq r$. A set of generators of the free part of $\OO^*$ is called a {\em system of
fundamental units} of $\OO$ and we write $\{\epsilon_1, \ldots, \epsilon_r\}$ for a given system of
fundamental units.

We now restrict to the case $r = 2$, i.e., unit rank~2.  Given any system of fundamental units
$\{\epsilon_1,\epsilon_2\}$, consider the $2 \times 2$ matrix $M = (v_i(\epsilon_j))_{1\leq i,j \leq
2}$.  If we transform $M$ into Hermite Normal Form, then the resulting matrix entries correspond to
valuations of another system of fundamental units, $\{\eta_1,\eta_2\}$. This system is independent
of the original system, and is unique up to constants in $\Fq^*$.  Furthermore, $\D_0^S =
\ideal{\infty_1-\infty_0, \infty_2-\infty_0}$ and $R =
|\D_0^S/\PP^S|=\det(M)=v_1(\eta_1)v_2(\eta_2)$.

For the remainder of this paper, we assume that $S$ contains an infinite place $\infty_0$ of degree
$1$, so that $f=1$ and $h = h_{\OO}R$.  In this case, $h_{\OO}$ is generally very small, so we
operate in another set of ideals called the {\em infrastructure} of $K$.  To that end, we require
the notion of a {\em distinguished} divisor and ideal.

\subsection{Distinguished Divisors and Infrastructure}
\label{S:distinguished}

Let $K$ be a cubic function field with an infinite place $\infty_0$ of degree $1$ and maximal order
$\OO$. A divisor $D$ of $K$ is said to be {\em finitely effective} if $D_S\geq 0$; that is,
$v_{\p}(D) \geq 0$ for all finite places $\p$ of $K$. Following \cite{b,gps,ericphd}, a finitely
effective divisor $D$ is defined to be {\em distinguished} if
\begin{enumerate}
  \item $D$ is of the form $D = D_S - \deg(D_S)\infty_0$, and
  \item if $E$ is any finitely effective divisor equivalent to $D$ with $\deg\left(E_S\right) \leq
  \deg\left(D_S\right)$ and $E^S \geq D^S$, then $D=E$.
\end{enumerate}
A fractional ideal $\f$ of $\OO$ is said to be \emph{distinguished} if $\Psi^{-1}(\f)$ is a
distinguished divisor.  Note that distinguished ideals are called {\em reduced} in
\cite{s00,ss00,s01,lsy,felix08,felixphd}.

A general treatment of infrastructures in function field extensions of arbitrary degree can be found
in \cite{felix08,felixphd}.  The cubic scenario was first presented in \cite{ss00,s01,lsy}, and we
use a description based on \cite{ericphd} here.

By \cite[Lemma 3.3.12 and Theorem 3.3.16]{ericphd}, if $K$ is a cubic function field with an
infinite place of degree $1$, then every divisor class contains at most one distinguished
divisor. (In fact, almost all divisor classes contain a distinguished divisor; see
\cite{felixholes}.)  This gives rise to the following definition. The (finite) set
$$\tilde{\R} := \left\{ \f\in\HH  \st \f \mbox{ is distinguished}  \right\}$$
is the {\em (principal) infrastructure} of $\OO$ (or of $K$). While we use an ideal-theoretic
definition of $\tilde{\R}$ here, the isomorphism $\Psi$ can be used to translate this into
divisor-theoretic language. In particular, $\tilde{\R}$ is in one-to-one correspondence with the set
of distinguished representatives of the kernel of the map $\J \to Cl(\OO)$.

Henceforth, we will restrict to the case $r=2$.  We consider the lattice $\Lambda :=
\ideal{(v_1(\eta_1),0),(v_1(\eta_2),v_2(\eta_2))} \subseteq \Z^2$. If $\f\in\tilde{\R}$, then there
is a function $\alpha \in K^*$ such that $\f = \ideal{\alpha^{-1}}$. The coset $(v_1(\alpha),
v_2(\alpha)) + \Lambda$ is uniquely determined by $\f$.  We define the {\em distance} of $\f$ to be
$\delta(\f) := (\delta_1(\f),\delta_2(\f)) + \Lambda := \left(v_1(\alpha),v_2(\alpha)\right) +
\Lambda$. Since $\delta : \tilde{\R} \to \Z^2/\Lambda$ is injective, $\tilde{\R}$ can be thought of
as a subset of $\Z^2/\Lambda$. In other words, $\tilde{\R}$ is structured as discrete points on the
surface of a torus.

In practice, we do not know $\Lambda$, and finding $\delta(\f)$ given only $\f$ is computationally
infeasible. We therefore define the \emph{(extended principal) infrastructure} as \[ \R := \{ (\f,
v) \in \tilde{\R} \times \Z^2 \mid \delta(\f) = v + \Lambda \}. \] For $\ai = (\f, v) \in \R$, call
$\delta(\ai) := (\delta_1(\ai), \delta_2(\ai)) := v$ the \emph{distance} and $\idealp(\ai) := \f$
the \emph{ideal part} of $\ai$. Finally, we will call $v_1(\eta_1) = \exp(\D_0^S/\PP^S)$ the {\em
exponent} of $\R$ and denote it $\exp(\R)$.  Since $R = v_1(\eta_1)v_2(\eta_2)$, we have
$\exp(\R)\mid R$.

\subsection{Infrastructure Arithmetic}
\label{S:inf-arith}

Infrastructures have two main operations: the {\em baby step} and {\em giant step} operations.
Roughly speaking, a baby step maps an infrastructure element to another element close to it, in
terms of distance, while a giant step reduces the product of two distinguished ideals.  We will also
describe a third operation called the {\em below} operation, which finds an infrastructure element
of (or close to and just below) a given distance.  Moreover, these operations can be computed
efficiently; for full details and proofs of this arithmetic in purely cubic function fields, we
refer the reader to \cite{ss00,s01,lsy,b,ericphd}.

In unit rank $2$ infrastructures, there are three types of baby steps as follows. Let $\f$ be a
distinguished ideal of $\OO$ and denote
\begin{eqnarray*}
\HH_i(\f) &=& \left\{\alpha \in \f \st v_i(\alpha) < 0, v_j(\alpha) \geq 0 \mbox{ for all } j\neq i, \right. \\ 
&&\left. \mbox{ and } v_j(\alpha) > 0 \mbox{ for at least one } j \neq i \right\} \nonumber\enspace .
\end{eqnarray*}
Following \cite{ss00,lsy}, let $i\in\{0,1,2\}$ and $\alpha, \beta \in K^*$.  Write $\alpha \geq_i
\beta$ if

$$\left(v_i(\alpha), v_{i+1}(\alpha), v_{i+2}(\alpha)\right) 
\geq_{lex} \left(v_i(\beta), v_{i+1}(\beta), v_{i+2}(\beta)\right) \enspace ,$$ where the subscripts
are considered modulo $3$.  The following theorem guarantees the existence and uniqueness (up to a
factor in $\Fq^*$) of the maximal element of $\HH_i(\f)$ under the ordering $\geq_i$.

\begin{theorem}[Theorem 3.4 of \cite{lsy}]
  \label{thm:H}
  Let $\OO$ be the maximal order of a purely cubic function field $K$ of unit rank 2 and $\f$ a
  distinguished fractional ideal of $\OO$.  For any $i\in\{0,1,2\}$, there exists an element $\phi =
  \phi_i(\f) \in \HH_i(\f)$, unique up to a factor in $\Fq^*$, such that $\phi \geq_i \alpha$ for
  all $\alpha \in \HH_i(\f)$. Furthermore, $\ideal{\phi^{-1}}\f$ is also a distinguished fractional
  ideal.
\end{theorem}
Let $i\in\{0,1,2\}$. If $\ai = (\f, v) \in\R$, $\phi = \phi_i(\f)$, and $\g = \ideal{\phi^{-1}}\f
\in \tilde{\R}$, then the operation $\ai \mapsto \bi := (\g, v + (v_1(\phi), v_2(\phi)))$ is called
a {\em baby step (in the $i$-direction)}, and we write $bs_i(\ai) = \bi$. See also
Figure~\ref{fig:babysteps} on how baby steps behave with high probability.

The giant step operation is analogous to multiplication. If $\ai_1 = (\f_1, v_1),\,\ai_2 = (\f_2,
v_2) \in\R$, then $\f_1\f_2$ is generally not distinguished. However, by \cite[Theorem
5.3.17]{ericphd}, there is a function $\psi\in K^*$ such that $v_i(\psi) \geq 0$, for each
$i=0,1,2$, and $v_0(\psi) + v_1(\psi) + v_2(\psi) \leq 2g$, yielding $\ai_1 * \ai_2 :=
(\ideal{\psi^{-1}}\f_1\f_2, v_1 + v_2 + (v_1(\psi), v_2(\psi))) \in \R$. Thus, $\delta(\ai_1*\ai_2)
= \delta(\ai_1) + \delta(\ai_2) + (v_1(\psi), v_2(\psi))$, so that $\delta(\ai_1*\ai_2) \gtrapprox
\delta(\ai_1)+\delta(\ai_2)$.  We call $*$ the {\em giant step} operation. Under $*$, $\R$ is an
abelian group-like structure, failing only associativity, and by \cite{felixholes}, existence of
inverses for very few elements.

A third required operation is the computation of the infrastructure element {\em below} any ordered
pair $(a,b)$ of integers $a,b\in\N$. This is the unique element $B(a,b):=\ai\in\R$ such that
$\delta(\ai) = (a - i, b)$ with $i \ge 0$ minimal.  From \cite{felixholes}, it follows that
$\delta(B(a,b)) = (a,b)$ with probability $1-O(1/q)$.

Navigating $\R$ is not as straightforward as the cyclic infrastructures of fields of unit rank
$1$. This is due to the existence of ``hidden'' elements and ``holes''.  An element in $\R$ is {\em
hidden} if it cannot be reached via baby steps. A \emph{hole} $d \in \Z^2$ is an element that does
not lie in the image of the distance map~$\delta$, i.e., there exists no element~$\ai \in \R$ with
$\delta(\ai) = d$. By \cite{felixholes}, the probability of encountering a hidden element or a hole
is $1 - O(1/q)$. Therefore, the distance advances effected by a baby step and a giant step are
generally predictable as follows.

\begin{proposition}
  \label{the-heuristic} Let $K/\Fq(x)$ be a cubic function field of genus $g$.
  \begin{enumerate}
    \item If $\ai\in\R$, then with probability $1-O(1/q)$, we have
    $$\delta(bs_i(\ai)) - \delta(\ai) = \left\{\begin{array}{ll}
      (1,0) & \mbox{ if } i = 0 \enspace , \\
      (-1,1) &  \mbox{ if } i=1 \enspace , \\
      (0,-1) &  \mbox{ if } i=2 \enspace .\\
    \end{array}\right.  $$
    \item If $\ai_1 = (\f,v_1),\,\ai_2 = (\f, v_2)\in\R$ and $\ai_1*\ai_2 =
    (\ideal{\psi^{-1}}\f_1\f_2, v_3)$, so that $\delta(\ai_1*\ai_2) = \delta(\ai_1) + \delta(\ai_2)
    + (v_1(\psi), v_2(\psi))$, then with probability $1-O(1/q)$, we have
    $$\ph(g):= v_1(\psi) = v_2(\psi) = \left\{\begin{array}{ll}
      \lfloor g/3\rfloor & \mbox{ if } g \not\equiv 2 \ppmod{3} \enspace , \\
      (g+1)/3 & \mbox{ if } g \equiv 2 \ppmod{3} \enspace .
    \end{array}\right.  $$
  \end{enumerate}
\end{proposition}
The proof follows from \cite{ericphd} using \cite{felixholes}. The first statement is visualized in
Figure~\ref{fig:babysteps}.

The majority of our computations will take place in the subset $\R_0 := \{\ai\in \R\st \delta_2(\ai)
= 0 \}$ of $\R$. As such, if we encounter an element $\bi\notin\R_0$, either via a baby step or a
giant step, then we must find an element in $\R_0$ close to $\bi$. Algorithm~\ref{alg:want2d0} finds
such an element $\ai\in\R_0$ with overwhelming probability and is based on Theorem~\ref{thm:H}. The
idea of the algorithm is to first find an element of non-negative $2$-distance. If the element has
positive $2$-distance at this point, then we expect that a step in the $0$-direction followed by a
series of steps in the $2$-direction produces an element $\ai \in \R_0$ with $\delta_1(\ai) \ge
\delta_1(\bi)$ and small $\delta_2(\ai) \ge 0$. If at this point, $\delta_2(\ai) \neq 0$, we repeat
the process until an element~$\ai \in \R_0$ is found. Since the number of holes is very small by
\cite{felixholes}, one iteration almost always suffices.

Figure~\ref{fig:want2d0} illustrates the most common scenario in which Algorithm~\ref{alg:want2d0}
is used, namely when we encounter a hole in $\R$ when taking a baby step in the $0$-direction.  This
baby step generally results in an element with a $2$-distance of $1$.

\tikzstyle{bstep} = [draw,thick,->]
\tikzstyle{bstepg} = [draw,thick,->,color=black!55]
\tikzstyle{inf} = [draw, shape = circle, fill=black, inner sep=0pt, minimum size = 0.2cm]
\tikzstyle{infg} = [draw, shape = circle, color=black!55, fill=black!55, inner sep=0pt, minimum size = 0.2cm]
\tikzstyle{hole} = [draw, shape = circle, fill=white, inner sep=0pt, minimum size = 0.2cm]

\begin{figure}[t]
\begin{center}
\caption{The Typical Baby Step Behaviors According to Proposition~\ref{the-heuristic}~(1)}
\label{fig:babysteps}
  \begin{tikzpicture}[scale=0.6, label distance=-1mm, node distance=0mm, >=stealth]
    \draw[->] (-0.75,-0.75) -- (-0.75,2.25);
    \draw[->] (-0.75,-0.75) -- (12.75,-0.75);
    \node (empty) at (12.75,-0.75) [label=right: \scriptsize$v_1$]{};
    \node (empty) at (-0.75,2.25) [label=above: \scriptsize$v_2$]{};
    
    \node[infg] (I00-0) at (0,0) {};
    \node[infg] (I10-0) at (1,0) {};
    \node[infg] (I20-0) at (2,0) {};
    \node[infg] (I01-0) at (0,1) {};
    \node[inf] (I11-0) at (1,1) [label=left: \scriptsize$\ai$]{};
    \node[inf] (I21-0) at (2,1) [label=right: \scriptsize$bs_0(\ai)$]{};
    \node[infg] (I02-0) at (0,2) {};
    \node[infg] (I12-0) at (1,2) {};
    \node[infg] (I22-0) at (2,2) {};
    \path[bstep] (I11-0) -- (I21-0) node[inner sep=1pt, pos=0.5, above] {\scriptsize $0$};

    \node[infg] (I00-1) at (5,0) {};
    \node[infg] (I10-1) at (6,0) {};
    \node[infg] (I20-1) at (7,0) {};
    \node[infg] (I01-1) at (5,1) {};
    \node[inf] (I11-1) at (6,1) [label=-45: \scriptsize$\ai$]{};
    \node[infg] (I21-1) at (7,1) {};
    \node[inf] (I02-1) at (5,2) [label=135: \scriptsize$bs_1(\ai)$]{};
    \node[infg] (I12-1) at (6,2) {};
    \node[infg] (I22-1) at (7,2) {};
    \path[bstep] (I11-1) -- (I02-1) node[inner sep=0pt, pos=0.5, below left] {\scriptsize $1$};

    \node[infg] (I00-2) at (10,0) {};
    \node[inf] (I10-2) at (11,0) [label=below: \scriptsize $bs_2(\ai)$]{};
    \node[infg] (I20-2) at (12,0) {};
    \node[infg] (I01-2) at (10,1) {};
    \node[inf] (I11-2) at (11,1) [label=above: \scriptsize $\ai$]{};
    \node[infg] (I21-2) at (12,1) {};
    \node[infg] (I02-2) at (10,2) {};
    \node[infg] (I12-2) at (11,2) {};
    \node[infg] (I22-2) at (12,2) {};
    \path[bstep] (I11-2) -- (I10-2) node[inner sep=1pt, pos=0.5, right] {\scriptsize $2$};
    
  \end{tikzpicture}
\\
\caption{Algorithm \ref{alg:want2d0} -- The Most Common Scenario}
\label{fig:want2d0}%
  \begin{tikzpicture}[scale=0.8, label distance=0mm, >=stealth]


    \draw[->] (-0.75,0) -- (-0.75,1.15);
    \draw[->] (-0.75,0) -- (4.85,0);
    \node (empty) at (4.75,0) [label=right: $v_1$]{};
    \node (empty) at (-0.75,1) [label=above: $v_2$]{};


    \node[infg] (I00) at (0,0) {};
    \node[infg] (I10) at (1,0) {};
    \node[infg] (I20) at (2,0) {};
    \node[hole] (I30) at (3,0) {};
    \node[inf] (I40) at (4,0) {};
    \node[infg] (I01) at (0,1) {};
    \node[infg] (I11) at (1,1) {};
    \node[hole] (I21) at (2,1) {};
    \node[inf] (I31) at (3,1) {};
    \node[inf] (I41) at (4,1) {};
    
    
    \node (dots) at (2.75,1.7) {\scriptsize Invoke Algorithm~\ref{alg:want2d0} here.};
    \draw[->] (2.75,1.5) -- (2.933,1.166);
    
    \path[bstepg] (I20) -- (I31) node[inner sep=-0.25pt, pos=0.5, above left] {\scriptsize$0$};
    \path[bstep] (I31) -- (I41) node[inner sep=1.25pt, pos=0.5, above] {\scriptsize$0$};
    \path[bstep] (I41) -- (I40) node[inner sep=0.5pt, pos=0.5, right] {\scriptsize$2$};
    \path[bstepg] (I00) to [out=45,in=135] (I10) node[inner sep=1.25pt, pos=0.5, above] {\scriptsize$0$};
    \path[bstepg] (I10) to [out=45,in=135] (I20) node[inner sep=1.25pt, pos=0.5, above] {\scriptsize$0$};
    
    
    \node[hole] at (7,1.5) {};
    \node[right] (dots) at (7.5,1.5) {\scriptsize hole};
    \node[inf] (Ileg) at (6.8,0.75) {};
    \node[infg] (Ileg) at (7.2,0.75) {};
    \node[right] (dots) at (7.5,0.75) {\scriptsize infrastructure elements};
    \path[bstep] (6.6,0) -- (7.4,0) node[inner sep=1.25pt, pos=0.5, above] {\scriptsize$i$};
    \node[right] (dots) at (7.5,0) {\scriptsize baby step in direction~$i$};

  \end{tikzpicture}
\end{center}
\end{figure}

\begin{algorithm}[t]
\caption{\textbf{($\mathbf{\textbf{red}_0}$)} Finding an Element in $\R_0$} \label{alg:want2d0}
\begin{algorithmic}[1]
\REQUIRE An element $\bi\in\R$ such that $\delta_2(\bi) \neq 0$.
\ENSURE An element $\reduce_0(\bi) := \ai \in\R_0$ close to $\bi$ such that $\delta_2(\ai) = 0$
\WHILEOL{$\delta_2(\bi)<0$}{$\ai := \bi$, $\bi := bs_1(\ai)$}
\WHILE[Now $\delta_2(\bi) \geq 0$.]{$\delta_2(\bi) > 0$} 
   \STATE $\ai := \bi$, $\bi := bs_0(\ai)$
   \WHILEOL{$\delta_2(\bi) > 0$}{$\ai := \bi$, $\bi := bs_2(\ai)$}
   \IFOL{$\delta_2(\bi) < 0$}{$\bi := \ai$}
\ENDWHILE
\RETURN $\ai := \bi$
\end{algorithmic}
\end{algorithm}

\section{The Kangaroo Method in $\R$}
\label{S:kangaroo}

If we are given integers $E,\,U\in\N$ such that the divisor class number $h \in (E-U,\,E+U)$, then
the Kangaroo method may be optimized to compute a multiple of $\exp(\R)$ with a probabilistic
running time of $O\bigl(\sqrt{U}\bigr)$ giant steps. While the Baby Step-Giant Step method is
generally faster than the Kangaroo method, the Kangaroo method is preferred for larger computations
because it requires very little storage and can be parallelized efficiently.  Specifically, we will
describe the parallelized Kangaroo method of van Oorschot and Wiener \cite{vow,st05} and explain
important improvements that apply in particular to operating in infrastructures of unit rank $2$.
After our description, we will optimize its running time in Theorem~\ref{thm:rooinf2}.  Later, we
will show how to determine the regulator $R$ and the divisor class number $h$ from $\exp(\R)$ in
many cases.

\begin{figure}[t]
\caption[]{The Kangaroo Method in Unit Rank $2$ Infrastructures}
\label{fig:roo2}
\begin{center}
  \begin{tikzpicture}[scale=0.31, label distance=0mm, >=stealth]
    
    \filldraw[black!40, dashed, fill=black!10] (29.17,0) -- (0,-5) -- (32,-5) -- (33.17,0);
    
    \filldraw[fill=black!20] (0,0) -- (8,0)  -- (10,4) -- (2,4) -- (0,0);
    \draw[line width=1.5pt] (0,0) -- (8,0);
    \draw[line width=1.5pt, dashed] (8,0) -- (9.25,0);
    \draw[dashed] (10,4) -- (11.25,4);
    \draw[line width=1.5pt, dashed] (10.75,0) -- (12,0);
    \draw[line width=1.5pt] (12,0) -- (20,0);
    \draw (20,0) -- (22,4) -- (14,4) -- (12,0);
    \draw[line width=1.5pt, dashed] (20,0) -- (21.25,0);
    \draw[line width=1.5pt, dashed] (22.75,0) -- (24,0);
    \draw[dashed] (12.75,4) -- (14,4);
    \draw[dashed] (22,4) -- (23.25,4);
    \draw[dashed] (24.75,4) -- (26,4);
    \draw[line width=5pt, color=black!20, cap=round] (29.4,0) -- (32.9,0);
    \draw[line width=1.5pt] (29.17,0) -- (24,0);
    \draw (24,0) -- (26,4) -- (34,4) -- (32,0);
    \draw[(-)] (29.17,0) -- (33.17,0);
    \draw[line width=1.5pt] (29.17,0) -- (33.17,0);
    
    \draw[->] (0,0) -- (0,3.5);
    \node (empty) at (0,3.5) [label=left: \scriptsize $v_2$]{};
    \draw[line width=1.5pt, ->] (33.17,0) -- (35,0);
    \node (empty) at (35,0) [label=above: \scriptsize $v_1$]{};
    
    \node (dots) at (5,2) {\scriptsize $\tilde{\R}$};
    \node (dots) at (10,0) {$\pmb{\cdots}$};
    \node (dots) at (12,4) {$\cdots$};
    \node (dots) at (22,0) {$\pmb{\cdots}$};
    \node (dots) at (24,4) {$\cdots$};
    \node (empty) at (0,0) [inner sep=0mm, label=below: \scriptsize $(0\mbox{, } 0)$]{};
    \node (empty) at (8,0) [inner sep=0mm, label=below: \scriptsize $(\exp(\R)\mbox{, } 0)$, label=above left: \scriptsize $\eta_1$]{};
    \node (empty) at (20,0) [inner sep=0mm, label=below: \scriptsize $(R\mbox{, } 0)$]{};
    \node (empty) at (32,0) [inner sep=0mm, label=below: \scriptsize $(h\mbox{, } 0)$, label=above left: \scriptsize $\epsilon$]{};
    
    \node (empty) at (2,4) [inner sep=0mm, label=above: \scriptsize $(v_1(\eta_2)\mbox{, } v_2(\eta_2))$, label=below right: \scriptsize $\eta_2$]{};
    \node (empty) at (15.5,0) [inner sep=0mm, label=below: \scriptsize $\R_0$]{};
    
    \fill (0,0) circle (0.2cm) (8,0) circle (0.2cm) (12,0) circle (0.2cm) (20,0) circle (0.2cm) (24,0) circle (0.2cm) (32,0) circle (0.2cm);
    \fill (2,4) circle (0.2cm) (10,4) circle (0.2cm) (14,4) circle (0.2cm) (22,4) circle (0.2cm) (26,4) circle (0.2cm) (34,4) circle (0.2cm);
    
    \draw[(-)] (0,-5) -- (32,-5);
    \draw[->] (16,-5) arc (180:0:1.666667);
    \draw[->] (19.333333,-5) arc (180:0:1);
    \draw[->] (21.333333,-5) arc (180:0:2);
    \draw[->] (25.333333,-5) arc (180:0:1.33333);
    \draw[->] (22.666667,-5) arc (180:0:1.333333);
    
    \fill[black!50] (16,-5) circle (0.2cm) (19.33333,-5) circle (0.2cm) (21.333333,-5) circle (0.2cm) 
    (25.3333333,-5) circle (0.2cm) (28,-5) circle (0.2cm);
    \fill (22.6666667,-5) circle (0.2cm);
    \fill[white] (28,-5) circle (0.2cm);
    \draw (28,-5) circle (0.2cm);
    \draw[->] (28,-6.6) -- (28,-5.4);
    
    \node (empty) at (28,-7) {\scriptsize trap};
    
    \node (empty) at (0, -5) [inner sep=0mm, label=below: \scriptsize $(E-U\mbox{, }0)$]{};
    \node (empty) at (16,-5) [inner sep=0mm, label=below: \scriptsize $(E\mbox{, }0)$]{};
    \node (empty) at (22.6666667,-5) [inner sep=0mm, label=below: \scriptsize $(h\mbox{, }0)$]{};
    \node (empty) at (32,-5) [inner sep=0mm, label=below: \scriptsize $(E+U\mbox{, }0)$]{};
    
  \end{tikzpicture}
\end{center}
\end{figure}
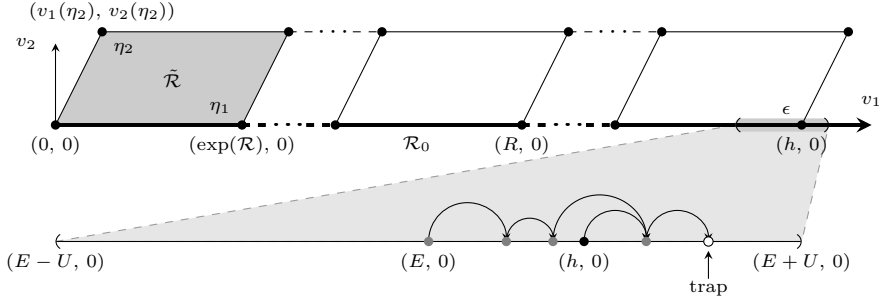

There are two key elements to adapting the Kangaroo algorithm to infrastructures of unit rank $2$
function fields. Firstly, the units of $\OO$ correspond to elements $(\OO, v) \in \R$, where $v \in
\Lambda$. Secondly, there exists a unit $\epsilon\in\OO^*$ such that $(v_1(\epsilon), v_2(\epsilon))
= (h,0)$; $\epsilon = \eta_1^i$, for some $i\in\N$.  Therefore, we restrict our search to
elements~$\ai$ with $\delta_2(\ai) = 0$, i.e., we operate in $\R_0 \subseteq \R$. In
Figure~\ref{fig:roo2}, we illustrate the Kangaroo method in our setting. The $v_1$ and $v_2$ axes
are labeled to give a reference for distance. The black dots correspond to units, with $\eta_1$,
$\eta_2$, and $\epsilon$ labeled. The infrastructure $\tilde{\R}$ is the gray parallelogram on the
left, with $\R_0$ the thick line at its base. Copies of $\tilde{\R}$ tile $\R$ in the
$v_1v_2$-plane.  A sample interval $(E-U,E+U)$ is shown in the top figure, highlighted in gray,
containing the unit $\epsilon$ at $(h,0)$. This interval is then expanded in the bottom figure to
show how the Kangaroo method proceeds.  Infrastructure elements (kangaroos) are initialized at
$(E,0)$ and $(h,0)$, which then jump, via baby steps and giant steps, along the $v_1$-axis until
their paths merge. The jumps are represented by the arcs.  Once these paths merge, a multiple of
$\exp(\R)$ can be determined, and from that we can determine $\exp(\R)$ itself. This often makes it
possible to determine $R$ and $h$ as well; in fact, in many cases $\exp(\R) = R = h$ (see the
discussion in Section~\ref{S:EU}). However, Figure~\ref{fig:roo2} illustrates the most general
situation.

We now describe in detail a modification of the parallelized Kangaroo method using notation similar
to that of \cite{st02,st05} for hyperelliptic function fields. Let $m$ be the (even) number of
available processors. The algorithm uses two {\em herds} of {\em kangaroos}: a herd
$\{T_1,\,\ldots,\,T_{m/2}\}$ of {\em tame} kangaroos and a herd $\{W_1,\,\ldots,\,W_{m/2}\}$ of {\em
wild} kangaroos.  A kangaroo is a sequence of elements in $\R$, and we write $T_j =
\{\tk_{A,j}\}_{A\in\N_0}$ and $W_k=\{\wk_{B,k}\}_{B\in\N_0}$, for $1\leq j,k\leq m/2$.  Each tame
and wild kangaroo is initialized via $\tk_{0,j} = B(E + (j-1)\nu,\, 0) \in \R_0$ and $\wk_{0,k} =
B((k-1)\nu,\, 0) \in \R_0$, respectively, for some small $\nu \in \Z$.  From these initial
positions, the kangaroos make {\em jumps} (i.e., baby and giant steps) in $\R_0$ until a {\em
collision} between a tame and a wild kangaroo occurs.  That is, the kangaroos jump until a tame and
a wild kangaroo have the same ideal part.  In this case, if $\idealp(\tk_{A,i}) =
\idealp(\wk_{B,i'})$, for some $A,B\in\N_0$ and $1\leq i,i'\leq m/2$, then $\delta_1(\tk_{A,i})
\equiv \delta_1(\wk_{B,i'}) \ppmod{\exp(\R)}$, so $h_0 := \delta_1(\tk_{A,i}) -
\delta_1(\wk_{B,i'})$ is a multiple of $\exp(\R)$.

To make the jumps, define a set of small (relative to $U$) random positive integers
$\{s_1,\,\ldots,\, s_{64}\}$, the {\em jump set} $J = \{\ai_1,\,\ldots,\,\ai_{64}\}$, where $\ai_i
=B(s_i-\ph(g),-\ph(g))$, for $1\leq i\leq 64$, and a hash function $w:\R \to \{1,\,\ldots,\,64\}$.
Also, for a real number $\tau\geq 1$, let $\SC_{\tau} \subseteq \R_0$ such that approximately every
$\tau$-th element of $\R_0$ belongs to $\SC_{\tau}$.  Each kangaroo jumps through $\R_0$ via an {\em
iteration} of a giant step and possibly one or more baby steps.  Initially, each tame and wild
kangaroo will take baby steps, if necessary, until it is in $\SC_{\tau}$. Then each kangaroo $\kk_l$
takes the giant step $\kk_{l+1} := \kk_l * \ai_{w(\kk_l)}$, for $l \geq 0$, followed by baby steps
in the $0$-direction, correcting via Algorithm~\ref{alg:want2d0} if necessary, until an element in
$\SC_{\tau}$ is found.

If there is a collision between two kangaroos of the same herd, then we must re-initialize one of
the two kangaroos.  If $\kk_l$ is one of two kangaroos in a collision, then choose a small $c \in
\N$, set $\kk_{l+1} := \kk_l* B(c,-\ph(g))$, and take baby steps until an element in $\SC_{\tau}$ is
found. Then $\kk$ continues jumping on its new path as usual. The other kangaroo may continue
without interruption.

Using the idea of van Oorschot and Wiener \cite{vow}, we will only store {\em distinguished points}
to reduce the storage requirement.  In order to avoid confusion in terminology, such points will be
called {\em (kangaroo) traps} instead. Let $\theta \in \N$ be a sufficiently large power of $2$ and
define another hash function $z:\R_0 \to \{0,\,\ldots,\,\theta-1\}$.  Set a trap, that is, store a
kangaroo $\kk$, if $z(\kk) = 0$.  Since kangaroos travel along the same path following a collision,
any collision will eventually land in a trap.

Finally, if there exist $a,\,b\in\N_0$, such that $b > 1$ and $h\equiv a\ppmod{b}$, then we make
adjustments to take advantage of this information.  First, change the estimate $E$ to $E - (E
\ppmod{b}) + a$, so that $E\equiv a\ppmod{b}$ for the revised value of $E$.  Next, choose $\nu$ and
the $s_i$ such that $b\mid \nu$ and $b\mid s_i$, for each $1\leq i\leq 64$.  Finally, restrict
$\SC_{\tau}$ to elements $\ai\in\R_0$ such that $\delta_1(\ai) \equiv a \ppmod{b}$ and require that
approximately every $b\tau$-th element of $\R_0$ lies in $\SC_\tau$.  The remaining initializations
and procedures are the same as before.

\begin{algorithm}[t]
\caption{Computing $h_0$ via the Kangaroo Algorithm} \label{alg:roo2}
\begin{algorithmic}[1]
\REQUIRE A purely cubic function field~$K/\F_q(x)$ of unit rank~$2$; $a,\,b \in \N_0$ such that $h \equiv a
\ppmod{b}$ (or $b = 1$ and $a = 0$ if no non-trivial $b$ is known); 
and an even integer $m$, the number of processors.%
\ENSURE A multiple~$h_0$ of $\exp(\R)$.%
\STATE Compute the genus $g$, choose $\rho$ from Table~\ref{tab:tau}, and choose $\ah := \ah(g)$ from Table~\ref{tab:alphas}.%
\STATE Set $\beta := \bigl[(m/2)\sqrt{(2\rho-1) \ah U}\bigr]-\rho+1$, $\nu := b\lfloor2\beta/(bm)\rfloor$,
$\theta := 2^{\left[\lg(\beta)/2\right]}$, $j := k := 0$.%
\STATE Choose random integers $g+1+\ph(g) \leq s_i\leq 2\beta$, with $1\leq i\leq 64$,
    such that $Mean\left(\{s_i\}\right) = \beta$ and $b\mid s_i$.%
\STATE Compute the jump set 
$J := \{\ai_1, \dots, \ai_{64} \}$, where $\ai_i := B(s_i-\ph(g),\,-\ph(g))$.%
\WHILE{$\delta_1(\ai_i) \neq s_i-\ph(g)$ for any $i$} %
  \STATE Replace $s_i:=s_i+b$ and recompute $\ai_i := B(s_i-\ph(g),\,-\ph(g))$ for inclusion in $J$.%
\ENDWHILE
\STATE Define hash functions $w:\R \to \{1,\,\ldots,\,64\}$ and $z: \R \to
\{0,\,\ldots,\,\theta-1\}$.%
\FOR{$i = 1$ to $m/2$}
   \STATE Initialize the tame kangaroos, $T_i$: $\tk_{0,i} := B(E+(i-1)\nu,\,0)$.%
   \STATE Initialize the wild kangaroos, $W_i$: $\wk_{0,i} := B((i-1)\nu,\,0)$.%
   \WHILEOL{$\tk_{0,i}\notin S_{\tau}$}{$\tk_{0,i} := \reduce_0(bs_0(\tk_{0,i}))$.}
   \WHILEOL{$\wk_{0,i}\notin S_{\tau}$}{$\wk_{0,i} := \reduce_0(bs_0(\wk_{0,i}))$.}
\ENDFOR
\WHILE{a collision between a tame and a wild kangaroo has not been found}
   \FOR{$i = 1$ to $m/2$}
      \IFOL{$z(\tk_{j,i}) = 0$ or $z(\wk_{k,i}) = 0$}{store the respective element(s).}
      \STATE Compute $\tk_{j+1,i} := \reduce_0\bigl(\tk_{j,i} * \ai_{w(\tk_{j,i})}\bigr)$ and 
             $\wk_{k+1,i} := \reduce_0\bigl(\wk_{k,i} * \ai_{w(\wk_{k,i})}\bigr)$. \label{step:wildgs}
      \WHILEOL{$\tk_{j+1,i}\notin S_{\tau}$}{$\tk_{j+1,i} := \reduce_0(bs_0(\tk_{j+1,i}))$.}
      \WHILEOL{$\wk_{k+1,i} \notin S_{\tau}$}{$\wk_{k+1,i} := \reduce_0(bs_0(\wk_{k+1,i}))$.}
      \STATE Increment $j:= j+1$ and $k:=k+1$.
   \ENDFOR
\ENDWHILE
\IFOL{$\tk_{A,\,i} = \wk_{B,\,i'}$}{\textbf{return}\ $h_0 := \delta_1(\tk_{A,i})-\delta_1(\wk_{B,i'})$.}
\end{algorithmic}
\end{algorithm}

In Algorithm~\ref{alg:roo2}, we formalize the procedures described above.  The following result is a
generalization of and an improvement upon similar ideas in \cite{st02,st05} and establishes optimal
choices for $\tau$ and the average jump distance $\beta = Mean(s_i)$ to minimize the expected
heuristic running time of the Kangaroo method. The proof is similar to the analogous result in cubic
function fields of unit rank $1$. We therefore omit the proof and refer the reader to
\cite{ericphd}.

\begin{theorem}
  \label{thm:rooinf2}
  Let $K/\Fq(x)$ be a purely cubic function field of unit rank $2$ such that $h \equiv a \ppmod{b}$
  for some $a,b\in\N$.  Then the expected heuristic running time, over all cubic function fields
  over $\Fq(x)$ of genus $g$, to compute a multiple $h_0$ of $\exp(\R)$ via Algorithm \ref{alg:roo2}
  is minimized by choosing $\tau = \rho/b$ and an average jump distance of $\beta =
  \bigl[(m/2)\sqrt{(2\rho-1) \alpha U}\bigr] - \rho+1$.  Here, $m$ is the (even) number of
  processors, $\rho = T_G/T_B$, $T_G$ and $T_B$ are the respective times required to compute a giant
  step and a baby step in $\R$, and $\alpha = \alpha(q,\,g) < 1/2$ is the mean value of $|h-E|/U$
  over all cubic function fields over $\Fq(x)$ of genus $g$. With these choices, the expected
  heuristic running time is 
  $\bigl(4\sqrt{\alpha U/(2\rho-1)} + \theta m + O(1)\bigr)(2-1/\rho)T_G$, as $q \to \infty$, where
  traps are set on average every $\theta$ iterations.
\end{theorem}

Following the recommendations given in \cite{teske,st05}, we make choices for the remaining
variables. First, we choose the $s_i$ randomly such that $g+1+\ph(g) \leq s_i \leq 2\beta$, for
$1\leq i\leq 64$. (The lower bound is an application of Theorem~5.3.10 of \cite{ericphd} to
guarantee that $(\OO,(0,0))\notin J$.)  The choice of $|J| = 64$ as a power of $2$ ensures that the
hash function $w$ is fast, and is large enough to obtain a sufficient level of randomization, but
small enough so that the space to store the jumps is not too large.  We also chose the spacing $\nu
\lessapprox 2\beta/m$.  For setting traps, we took $\theta = 2^{[\lg(\beta)/2]+c} =
O\bigl(\sqrt[4]{U}\bigr)$, for some small integer $c$.  The hash functions~$w$ and $z$ are defined
using an $\F_q[x]$-module representation of the ideal component of each element in $\R$; for
details, see \cite{ericphd}.

Table~\ref{tab:tau} lists values of $\rho$ for various unit rank $2$ situations of genera $2\leq
g\leq 7$. In each case, we computed the ratios using $10^6$ baby steps and $10^6$ giant steps in a
function field $\Fq(x,y)$ with $q = 10^8 + 39$ and $y^3 = GH^2$, where $G$ and $H$ were random,
monic, co-prime, irreducible polynomials with $\deg(G) \geq \deg(H)$.

In the next section, we briefly review the method of \cite{ss07} implemented here to compute the
divisor class number of a cubic function field.

\newsavebox{\tmpbox}
\savebox{\tmpbox}{\begin{minipage}[t]{9cm}\small Giant Step to Baby Step Ratio $\rho = T_G/T_B$ and
Estimate $\ah(g)$ of $Mean(|h-E|/U)$ \end{minipage}}
\begin{table}[t]
\begin{center}%
\begin{tabular}{|c|c|c||c|||c|c||c|c|c|r|r||c|}
\cline{1-7} \cline{9-12}
$g$ & $\deg(G)$ & $\deg(H)$ & $\rho$ & $\deg(G)$ & $\deg(H)$ & $\rho$ & & $g$ & $q$ & $\lambda$ & $\ah(g)$\\
\cline{1-7} \cline{9-12}
2 & 2 & 2 & 3.04839 & & & & \phantom{xx} & & & & \\
\cline{1-7} \cline{9-12}
3 & 4 & 1 & 3.02410 & & &             & & 3 & 100003 & 1 & 0.27187490 \\
\cline{1-7} \cline{9-12}
4 & 6 & 0 & 4.03846 & 3 & 3 & 4.42018 & & 4 & 10009  & 1 & 0.19186318 \\
\cline{1-7} \cline{9-12}
5 & 5 & 2 & 5.61416 & & &             & & 5 & 997    & 2 & 0.19190607 \\
\cline{1-7} \cline{9-12}
6 & 7 & 1 & 5.96440 & 4 & 4 & 6.38660 & & 6 & 463    & 2 & 0.15975657 \\
\cline{1-7} \cline{9-12}
7 & 9 & 0 & 7.87264 & 6 & 3 & 8.21655 & & 7 & 97     & 2 & 0.12602172 \\
\cline{1-7} \cline{9-12}
\end{tabular}%
\end{center}%
\caption{\usebox{\tmpbox}}%
\label{tab:tau}
\label{tab:alphas}
\end{table}

\section{Computing $h$ and $R$ -- the Idea}
\label{S:EU}

Algorithm~\ref{alg:h-idea} lists the three main phases of the method of Scheidler and Stein
\cite{ss07,ss08} to compute a multiple $h_0$ of $\exp(\R)$. If $\exp(\R)$ is large enough, then
these three steps determine the divisor class number $h$ of a cubic function field.  Step~4
determines $\exp(\R)$ and in certain cases, Step~5 computes the regulator~$R$ and ideal class number
$h_{\OO}$ of $\OO$.

If $\exp(\R) \le 2 U$, then there may be more than one multiple of $\exp(\R)$ in the interval~$(E -
U, E + U)$, in which case $h$ cannot be determined.  Nonetheless, $h$ is limited to a smaller
subset, since $\exp(\R) \mid h$. We know that $R/\exp(\R)$ is a divisor of $d:=\gcd(\exp(\R),
h/\exp(\R))$; if $d=1$, then $R = \exp(\R)$.

By results of Achter and Pries \cite{achter-districlassgroups,achterpries-intmonodromy}, the class
numbers of purely cubic function fields of genus~$g$ over a finite field~$\F_q$ behave like random
integers in the Hasse-Weil interval $[(\sqrt{q} - 1)^{2 g}, (\sqrt{q} + 1)^{2 g}]$ with respect to
divisibility. Therefore, the class numbers are very often square-free, whence the divisor class
group~$\J$ is cyclic. In that case, $\exp(\R) = R$, and if one assumes that $\D_0^S/\PP^S$ behaves
like a random subgroup of $\J$, then $R$ is large. Therefore, Algorithm~\ref{alg:h-idea} can
determine $R$ and $h$ in many cases.

For details on how to compute $E$ and $U$ in Step 1, along with the complete analysis of the running
time of Algorithm~\ref{alg:h-idea}, see \cite{ss07,ss08}. Further implementation details may be
found in \cite{ericphd}.  Here we merely state that by \cite{ss07,ss08}, for $g \geq 3$, the
complexity of Step~1 of Algorithm~\ref{alg:h-idea} is $O\left(q^{[(2g-1)/5]+\varepsilon(g)}\right)$
giant steps, as $q \to \infty$, where $-1/4 \le \varepsilon(g) \le 1/2$. If $g \leq 2$, then there
is no asymptotic improvement in using the bounds described in \cite{ss07,ss08} versus the Hasse-Weil
bounds.

Next, we discuss some practical issues arising in our implementation of Algorithm~\ref{alg:h-idea}.
We omit details on Step 3 since they were already given in Section~\ref{S:kangaroo}.

\section{Implementation Details}
\label{S:details}

\subsection{Implementation Details for Phase 2}
\label{S:implementation2}

\begin{algorithm}[t]
\caption{Computing $h$ and/or $R$ -- the Idea} \label{alg:h-idea}
\begin{algorithmic}[1]
\STATE Determine $E,\,U\in\N$ such that $h \in (E-U,\,E+U)$.%
\STATE Determine extra information about $h$ such as congruences or the distribution of 
$h$ in the interval $(E-U,\,E+U)$.%
\STATE Compute a multiple $h_0$ of $\exp(\R)$ via Algorithm~\ref{alg:roo2}.%
\STATE Compute $R^*:=\exp(\R)$ via Algorithm \ref{alg:S-regulator2}.
\STATE If $R^* > 2 U$, let $h$ be the unique multiple of $R^*$ in $(E-U,\,E+U)$.
If $\gcd(R^*, h/R^*) = 1$, then $R = R^*$ and $h_{\OO} = h/R^*$.
\end{algorithmic}
\end{algorithm}

For Phase 2 of Algorithm~\ref{alg:h-idea}, we use extra information about $h$ to effectively reduce
the size of the interval, $(E-U,\,E+U)$, determined in Phase 1.  The method to compute $E$ and $U$
uses a truncated Euler product representation of the zeta function of the function field, and we
consider finite places (i.e., monic irreducible polynomials) up to a degree bound $\lambda$.  It has
been shown in both the quadratic and cubic function field cases that $h$ is not uniformly
distributed in this interval, and tends to be close to the approximation $E$ \cite{st02b,ericphd}.

Let $\alpha(q,\,g) = Mean(|h-E|/U)$, where the mean is taken over all cubic function fields of genus
$g$ over $\Fq(x)$.  In Theorem~\ref{thm:rooinf2}, we described how to apply $\alpha(q,g)$ to
minimize the expected running time of Algorithm~\ref{alg:roo2}.  For a fixed genus $g$, we assume
that the limit $\alpha(g) = \lim_{q\to\infty}\alpha(q,\,g)$ exists, as is the case for hyperelliptic
function fields \cite{st02b}.  However, $\alpha(q,\,g)$ and $\alpha(g)$ are very difficult to
compute precisely, so instead we applied approximations $\ah(g)$ of $\alpha(g)$ for $3 \leq g \leq
7$.  Table~\ref{tab:alphas} (Table 6.5 of \cite{ericphd}) lists these approximations for selected
values of $g$, based on a sampling of $10000$ cubic function fields of genus $g$ over a fixed field
$\Fq$.  However, these averages may be applied to cubic function fields over any finite field.  In
Table~\ref{tab:alphas}, $\lambda$ is the degree bound used to compute the estimate $E$.

A second component of Phase 2 of Algorithm~\ref{alg:h-idea} finds information about $h$ modulo small
primes. In \cite{btw}, Bauer, Teske, and Weng consider purely cubic function fields defined by
Picard curves. In this case, they proved the following result about $h$ modulo powers of $3$.

\begin{proposition}[Lemma 2.2 of \cite{btw}]
  \label{lem:3powerdivh}
  Let $K = \Fq(x,y)$ be the function field of a Picard curve $C:y^3 = F(x)$, where $q\equiv
  1\ppmod{3}$.  If $F$ has $k$ distinct irreducible factors over $\Fq[x]$, then $3^{k-1}\mid h$. If
  $F$ is irreducible, then $h \equiv 1 \ppmod{3}$.
\end{proposition}

The genus $3$ curves we used in our computations are birationally equivalent to Picard curves, so we
applied this proposition to these curves.

\subsection{Implementation Details for Phase 4}
\label{S:implementation4}

\begin{algorithm}[t]
\caption{Computing $\exp(\R)$: Step 4 of Algorithm~\ref{alg:h-idea}} 
\label{alg:S-regulator2}
\begin{algorithmic}[1]
\REQUIRE A multiple $h_0$ of $\exp(\R)$ and a lower bound $l$ of $\exp(\R)$.
\ENSURE The exponent $\exp(\R)$ of the infrastructure $\R$.
\STATE Set $h^* := 1$.%
\STATE Factor $h_0 = \prod_{i=1}^k p_i^{a_i}$. \label{step:factor}%
\FOR {$i=1$ to $k$}
   \IF {$p_i < h_0/l$}
      \STATE Find $1\leq e_i \leq a_i$ minimal such that $\idealp(B\left(h_0/p_i^{e_i},0\right)) \neq \OO$. \label{step:below}
      \STATE Set $h^* := p_i^{e_i-1}h^*$.
   \ENDIF
\ENDFOR
\RETURN $\exp(\R) = R^* := h_0/h^*$.
\end{algorithmic}
\end{algorithm}

Algorithm~\ref{alg:S-regulator2} outlines the procedure for Step 4 of Algorithm~\ref{alg:h-idea}.
This step will determine $\exp(\R)$ given a multiple $h_0$ of $\exp(\R)$.  Here, we adapt Algorithm
4.4 of \cite{sw} to the case of cubic function fields of unit rank $2$, using the fact that
$\exp(\R)$ 
is the smallest factor $R^*$ of $h_0$ such that $\idealp(B(R^*,0)) = \OO$.  Recall that $B(a,0)$,
for $a\in\N$, may be impossible to determine because of a hidden element having distance $(a,0)$.
Nevertheless, the probability of this occurring is negligible, so that we can assume that Algorithm
\ref{alg:S-regulator2} produces the correct output.

We briefly comment on the running time of Algorithm~\ref{alg:S-regulator2} relative to the running
time of Algorithm~\ref{alg:h-idea}, especially in light of the factorization in
Step~\ref{step:factor}. First, current heuristic methods to factor the integer $h_0$ require a
subexponential number of bit operations in $\log(h_0)$.  Furthermore, Step 3 only requires a
polynomial number (in $g$ and $\log(q)$) of infrastructure operations. Therefore, determining
$\exp(\R)$ from $h_0$ will not dominate the overall running time of Algorithm~\ref{alg:h-idea}.  The
class numbers that we found required only a few seconds to factor. In fact, we simply used a basic
implementation of Pollard's Rho method for factoring~\cite{pollardrho}.

\section{Computational Results}
\label{S:results}

In this section, we tested the practical effectiveness of the Kangaroo algorithm to compute the
divisor class number and extracted the ideal class number and regulator of six purely cubic function
fields of unit rank $2$: five of genus $3$ and one of genus $4$. We remark that this is the first
time that Algorithm~\ref{alg:h-idea} has been implemented for cubic function fields of unit rank
$2$.

The genus $3$ curves that we used for the examples in this section were each of the form $C_i:y^3 =
G_i(x)x^2$, where
$$\begin{array}{ll}
G_1(x) = x^4 + 858028x^3 + 786068x^2 + 69746x + 675670 \enspace , \\
G_2(x) = x^4 + 9655935x^3 + 8633555x^2 + 1319425x + 1437614 \enspace , \\
G_3(x) = x^4 + 63268943x^3 + 53257730x^2 + 59385220x +  16188628 \enspace , \\
G_4(x) = x^4 + 834364201x^3 + 8363484x^2 + 953863416x + 850202733 \enspace , \\
G_5(x) = x^4 + 9994854268 x^3 + 7631258748 x^2 + 7469686108 x + 292775976 \enspace ,
\end{array}$$
and the genus $4$ curve was of the form $C_6: y^3 = G_6(x)$, where
$$\begin{array}{ll}
G_6(x) = x^6 + 4207x^5 + 3340x^4 + 9858x^3 + 7507x^2 + 36x + 1019 \enspace .
\end{array}$$
Each $G_i$ is irreducible over the field $\Fq$ used, and $q\equiv 1\ppmod{3}$ is prime. 

In Table~\ref{tab:r2results}, we list the ideal class number $h_{\OO}$, the regulator $R$, and the
ratio $|h-E|/U$ for these six examples.  For the genus $3$ examples, we used $\rho = 3.02410$ and
$\tau = \rho/3 = 1.00803$, and for the genus $4$ example, we used $\rho = \tau = 4.03846$.  Based on
the last column, we see that the estimate $E$ was better than average except for the computations
with curves $C_3$ and $C_4$.  The $C_4$ through $C_6$ examples were computed via a parallelized
approach, using up to $64$ processors.  The largest divisor class number we computed had $31$
decimal digits.

Data from the Kangaroo computations is given in Table~\ref{tab:r2roodata}. Here, ``BS Jumps'' and
``GS Jumps'' refer to the respective number of baby steps and giant steps computed using the
Kangaroo method in each example, $\lg\theta$ is the base $2$ logarithm of the value of $\theta$ used
for setting traps, ``Traps'' is the total number of traps that were set, $m$ is the number of
processors (or kangaroos) that were used, ``Coll.'' is the number of useless collisions in the given
example, and ``Time'' refers to the total time taken by the computation in minutes, hours, and days.
For timing and technical considerations, we implemented our algorithms in C++ using NTL, written by
Shoup \cite{ntl}, compiled using {\tt g++}, and run on IBM cluster nodes with Intel Pentium 4 Xeon
$2.4$ GHz processors and $2$ GB of RAM running Redhat Enterprise Linux 3.

\begin{table}[H]
$$\begin{array}{|l|r|l||r|r||c|}
\hline
\mbox{Curve} & q & g & h_{\OO} & R & |h-E|/U \\
\hline
C_1 & 1000003 & 3 & 1 &           1002847489604613721  & 0.1498574 \\
C_2 & 10000141 & 3 & 1 &        1000397435760158462929  & 0.1263140 \\
C_3 & 100000039 & 3 & 1 &     1000094985874807321192993  & 0.3799612 \\
C_4 & 1000000009 & 3 & 1 &  1000036037504733195527721763  & 0.3814163 \\
C_5 & 10000200031 & 3 & 1 & 1000028959108091361595659615907 &  0.2262216 \\
\hline
C_6 & 10009 & 4 & 1 &            10081785007075827  & 0.1218925 \\
\hline
\end{array}$$
\caption{Regulators and Ideal Class Numbers}
\label{tab:r2results}

$$\begin{array}{|l|r|r||r|r||r|r||r|r|r|}
\hline
\mbox{Curve} & q & g & \mbox{BS Jumps} & \mbox{GS Jumps} & \lg\theta & \mbox{Traps} & m & \mbox{Coll.} & \mbox{Time}\\
\hline
C_1 & 1000003 & 3 & 18825 &  2353928 & 10 & 2337 &  2 & - &  98.5\, m \\
C_2 & 10000141 & 3 & 72149 &  9040560 & 12 & 2198 &  2 & - &  5.67\, h\\
C_3 & 100000039 & 3 &  1537984 & 192895918 & 14 & 11768 &  2 & - &  6.84\, d \\ 
C_4 & 1000000009 & 3 &  12942441 & 1624509536 & 18 & 6147 & 40 & 0 &  37.4\, d \\
C_5 & 10000200031 & 3 & 404518765 & 50757901157 & 20 & 46235 & 64 & 7 & 1543 \, d \\
\hline
C_6 & 10009 & 4 & 1814587 & 596886 & 12 & 109 &  6 & 2 &  80.3\, m \\
\hline
\end{array}$$
\caption{Regulator Computation Data}
\label{tab:r2roodata}
\end{table}

\section{Conclusions and Future Work}
\label{S:conclusions}

Using current implementations of the arithmetic in the infrastructure of a purely cubic function
field of unit rank $2$, divisor class numbers and regulators up to $31$ digits were computed using
the method of Scheidler and Stein \cite{ss07} and the Kangaroo algorithm as a subroutine. The
largest example among these class numbers and regulators was the largest ever computed for a
function field of genus at least $3$, with the exception of function fields defined by a Picard
curve. This was also the first time that a ``square-root'' algorithm was efficiently applied to the
infrastructure $\R$ of a global field of unit rank $2$.  Moreover, we made improvements to the
Kangaroo method in $\R$ by showing how to take advantage of information on the congruence class of
the divisor class number and how to use the ratio $\rho = T_G/T_B$ more effectively.

A procedure to determine the regulator given the divisor class number and infrastructure exponent
using methods of Buchmann, Jacobson, and Teske \cite{bjt,teske98} is work in progress.  In addition,
efficient ideal and infrastructure arithmetic needs to be developed for arbitrary (i.e., not
necessarily purely) cubic function fields as well as for characteristic $2$ and $3$ in order to
apply this method to such function fields.  Finally, it is unknown if we can take advantage of the
torus structure of $\R$ to compute $R$ using more efficient techniques.

\end{document}